\documentclass[12pt,a4paper,draft]{article}
\usepackage[english,russian]{babel}
\usepackage{amssymb, amsmath, latexsym, amsfonts, amsthm}
\begin{document}
\renewcommand{\refname}{References}
\newtheorem{theorem}{Theorem}
\newtheorem{lemma}{Lemma}
\newtheorem{corollary}{Corollary}
\newtheorem{proposition}{Proposition}
\begin{center}
On a space of entire functions and its Fourier transformation
\end{center}
\begin{center}
I Kh. Musin\footnote {E-mail: musin@matem.anrb.ru}
\end{center}
\begin{center}
Institute of Mathematics with Computer Centre of Ufa Scientific Centre of Russian Academy of Sciences, 

Chernyshevsky str., 112, Ufa, 450077, Russia
\end{center}

\vspace {0.2cm}

\renewcommand{\abstractname}{}
\begin{abstract}
{\sc Abstract}. A space of entire functions of several complex variables rapidly decreasing on ${\mathbb R}^n$ and such that their growth along $i{\mathbb R}^n$ is majorized with a help of a family of weight functions (not radial in general) is considered in the paper. For this space an equivalent description in terms of estimates on all partial derivatives of functions on ${\mathbb R}^n$ and Paley-Wiener type theorem are given. 

\vspace {0.2cm}






\vspace {0.3cm}

\end{abstract}

\vspace {1cm}

{\bf 1. On the problem}. In the 1950's the study of $W$-type spaces started  with the works of B.L. Gurevich \cite {Gur1}, \cite {Gur2} and I.M. Gelfand and G.E. Shilov \cite {GS1}, \cite {GS2}. They described them by means of the Fourier transformation and then applied this description to study the uniqueness of the Cauchy problem of partial differential equations and their systems. 
These spaces generalize Gelfand-Shilov spaces of $S$-type  \cite {GS1}. So they are often called  Gelfand-Shilov spaces of $W$-type. 

Recall the definition of the Gelfand-Shilov spaces of $W$-type. 
Let $M$ and $\Omega$ be differentiable functions on $[0, \infty)$ such that 
$M(0)= \Omega(0) = M'(0) = \Omega'(0) = 0$ and their derivatives are continuous, increasing and tending to infinity. Then:

$W_M$ is a space of infnitely differentiable functions
$f$ on ${\mathbb R}^n$ satisfying the estimate 
$$
\vert (D^{\beta} f) (x) \vert \le C_{\beta} e^{-M(a\Vert x \Vert)}, \ x \in {\mathbb R}^n,
$$
where $a$ is some positive constant;

$W^{\Omega}$ is a space of entire functions $f$ on ${\mathbb C}^n$ satisfying the estimate 
$$
\vert \zeta^{\alpha} f(\zeta) \vert \le C_{\alpha} 
e^{\Omega(b \Vert \eta \Vert)}, \ \zeta = \xi +i \eta \in {\mathbb R}^n + i{\mathbb R}^n,
$$ 
for some $b > 0$;

$W_M^{\Omega}$ is a space of entire functions $f$ on ${\mathbb C}^n$ satisfying
the estimate 
$$
\vert f(\xi +i \eta) \vert \le C e^{-M(\Vert \xi \Vert) + \Omega(b \Vert \eta \Vert)}, \ \xi, \eta \in {\mathbb R}^n,
$$ 
for some positive $a, b$ and $C$.

$W$-type spaces and some their generalizations were studied by many mathematicians. For example, new characterizations of $W$-type spaces and their generalizations have been given  by J. Chung, S.Y. Chung and D. Kim \cite {C-C-K 1}, \cite {C-C-K 2},
R.S. Pathak and S.K. Upadhyay \cite {P-U}, S.K. Upadhyay \cite {U} 
(in terms of Fourier transformation), 
N.G. De Bruijn \cite {Br}, A.J.E.M. Janssen and S.J.L. van Eijndhoven \cite {J-E}, Jonggyu Cho \cite {Cho}
(by using the growth of their Wigner distributions). R.S. Pathak \cite {P} and S.J.L. van Eijndhoven and M.J. Kerkhof \cite {E-K} introduced new spaces of W-type and investigated the behaviour of the Hankel transformation over
them (see also  \cite {B-R}). New $W$-type spaces introduced by V.Ya. Fainberg, M.A. Soloviev \cite {FS} turned out to be useful for nonlocal theory of highly singular quantum fields. $W^{\Omega}$-type spaces generated by weight functions (radial, not convex in general) satisfying some natural conditions were studied in  \cite {M-M 1}, \cite {M-M 2}.

Here we consider $W^{\Omega}$-type spaces of entire functions constructed with a help of some family of weight functions of more general nature than in \cite {M-M 1}, \cite {M-M 2} (they are not convex, not radial in general) and announce our recent results devoted to such spaces. 
Full proofs of results are given in a paper submitted to the journal "`Concrete operators"'. 

We work with the space $E(\varPhi)$ of entire functions that is defined as follows.
Let $n \in {\mathbb N}$, $H({\mathbb C}^n)$ be the space of entire functions on ${\mathbb C}^n$,  
$\Vert u \Vert$ be the Euclidean norm of $u \in {\mathbb R}^n ({\mathbb C}^n)$. 
Denote by ${\mathcal A}$ the set of all real-valued functions $g \in C({\mathbb R}^n)$ satisfying the following three conditions:

1). $g(x) = g(\vert x_1 \vert, \ldots , \vert x_n \vert), \ x = (x_1, \ldots , x_n)\in {\mathbb R}^n$;

2). restriction of $g$ to $[0, \infty)^n$ is nondecreasing in each variable;

3). $\displaystyle \lim_{x \to \infty} \frac {g(x)}{\Vert x \Vert}= + \infty$.

\noindent Let $\varPhi =\{\varphi_m\}_{m=1}^{\infty}$ be a subset of ${\mathcal A}$ consisting of 
functions $\varphi_m$ satisfying the following condition:

$i_0$). for each $m \in {\mathbb N}$ and each $A > 0$ there exists a constant $C_{m, A} > 0$ such that 
$$
\varphi_m(x) + A \ln (1 + \Vert x \Vert) \le \varphi_{m+1}(x) + C_{m, A}, \ x \in {\mathbb R}^n.
$$
\noindent For each $\nu \in {\mathbb N}$ and $k \in {\mathbb Z}_+$ define the space
$$
E_k(\varphi_{\nu}) = \{f \in H({\mathbb C}^n): p_{\nu, k}(f) = \sup_{z \in {\mathbb C}^n} 
\frac 
{\vert f(z)\vert (1 + \Vert z \Vert)^k}
{e^{\varphi_{\nu} (\vert Im \ z_1 \vert, \cdots, \vert Im \ z_n \vert)}} < \infty \}.
$$
Let $E(\varphi_{\nu})= \bigcap \limits_{k=0}^{\infty} E_k(\varphi_{\nu})$, $E(\varPhi)= \bigcup \limits_{\nu=1}^{\infty} E(\varphi_{\nu})$. 
With usual operations of addition and multiplication by complex numbers 
$E(\varphi_{\nu})$ and $E(\varPhi)$ are linear spaces. 
Since $p_{\nu, k}(f) \le p_{\nu, k+1}(f)$ for $f \in E_{k+1}(\varphi_{\nu})$ then $E_{k+1}(\varphi_{\nu})$ is continuously embedded in $E_k(\varphi_{\nu})$. 
Endow $E(\varphi_{\nu})$ with a projective limit topology of spaces $E_k(\varphi_{\nu})$. 
Note that if $f \in E(\varphi_{\nu})$ then 
$p_{\nu+1, k}(f) \le e^{C_{\nu, 1}} p_{\nu, k}(f)$ for each $k \in {\mathbb Z}_+$. 
Hence, $E(\varphi_{\nu})$ is continuously embedded in $E(\varphi_{\nu + 1})$ for each $\nu \in {\mathbb N}$. 
Supply $E(\varPhi)$ with a topology of an inductive limit of spaces $E(\varphi_{\nu})$.

In the paper description of the space $E(\varPhi)$ in terms of estimates 
on partial derivatives of functions on ${\mathbb R}^n$ and description of Fourier transforms of functions of $E(\varPhi)$ under additional conditions on $\varPhi$ is given. Results of the paper could be useful in harmonic analysis, theory of entire functions, convex analysis, in the study of partial differential and pseudo-differential operators.

{\bf 2. Some notations 	and definitions}.  
For $u=(u_1, \ldots , u_n) \in {\mathbb R}^n \ ({\mathbb C}^n)$, $v=(v_1, \ldots , v_n) \in {\mathbb R}^n \ ({\mathbb C}^n)$ \ 
$\langle u, v \rangle  = u_1 v_1 + \cdots + u_n v_n$. 

For $\alpha = (\alpha_1, \ldots , \alpha_n) \in {\mathbb Z}_+^n$, 
$x =(x_1, \ldots , x_n) \in {\mathbb R}^n$, 
$z =(z_1, \ldots , z_n) \in {\mathbb C}^n$ \ $\vert \alpha \vert = \alpha_1 + \ldots  + \alpha_n$, 
$\alpha! = \alpha_1! \cdots \alpha_n!$, 
$x^{\alpha} = x_1^{\alpha_1} \cdots x_n^{\alpha_n}$, 
$D^{\alpha}=
\frac {{\partial}^{\vert \alpha \vert}}{\partial x_1^{\alpha_1} \cdots \partial x_n^{\alpha_n}}$ .

$\Pi_n:=[0, \infty)^n$.

If $\Pi_n \subseteq X \subset {\mathbb R}^n$ then for a function $u$ on $X$ denote by $u[e]$ the function defined by the rule:
$u[e](x) = u(e^{x_1}, \ldots, e^{x_n}), \ x = (x_1, \ldots , x_n) \in  {\mathbb R}^n$, and by $R(u)$ denote the restriction of $u$ to $\Pi_n$.
For brevity denote $\varphi_m [e]$ by $\psi_m$.

The Young-Fenchel conjugate of a function $g:{\mathbb R}^n \to [-\infty, + \infty]$ is the function 
$g^*:{\mathbb R}^n \to [-\infty, + \infty]$ defined by
$
g^*(x) = \displaystyle \sup \limits_{y \in {\mathbb R}^n}(\langle x, y \rangle - g(y)), \ x \in {\mathbb R}^n. 
$

{\bf 3. Main results}. Let $\Psi^*=\{\psi_{\nu}^*\}_{\nu=1}^{\infty}$. For each $\nu \in {\mathbb N}$ and $m \in {\mathbb Z}_+$ let 
$$
{\mathcal E}_m(\psi_{\nu}^*) =\{f \in  C^{\infty}({\mathbb R}^n): 
{\cal R}_{m, \nu}(f)= \sup_{x \in {\mathbb R}^n, \alpha \in {\mathbb Z}_+^n} 
\frac {(1+ \Vert x \Vert)^m \vert (D^{\alpha}f)(x) \vert}{\alpha! e^{-\psi_{\nu}^*(\alpha)}} < \infty \}.
$$
Let ${\mathcal E}(\psi_{\nu}^*) = \bigcap \limits_{m=0}^{\infty}{\mathcal E}_m(\psi_{\nu}^*)$, 
${\mathcal E}(\Psi^*) = \bigcup \limits_{\nu=1}^{\infty}{\mathcal E}(\psi_{\nu}^*)$.

The first two theorems are aimed to characterize functions of the space $E(\varPhi)$ 
in terms of estimates of their partial derivatives on ${\mathbb R}^n$.   

\begin{theorem} 
For each $f \in  E(\varPhi)$ its restriction $f_{|{\mathbb R}^n}$ 
to ${\mathbb R}^n$ belongs to ${\mathcal E}(\Psi^*)$.
\end{theorem}

\begin{theorem} 
Let the family $\varPhi$ satisfies the additional conditions:

$i_1)$. for each $m \in {\mathbb N}$ there exist constants $\sigma_m > 1$ and $\gamma_m > 0$ such that
$$
\varphi_m(\sigma_m x) \le \varphi_{m+1}(x) + \gamma_m, \ x \in {\mathbb R}^n;
$$

$i_2)$. for each $m \in {\mathbb N}$ there exists a constant $K_m > 0$ such that 
$$
\varphi_m(x + \xi) \le \varphi_{m+1}(x) + K_m, \ x \in \Pi_n, \xi \in [0, 1]^n.
$$

Then each function $f \in {\mathcal E}(\Psi^*)$ admits (the unique) extension to entire function belonging to $E(\varPhi)$.
\end{theorem}

Further, for $\nu \in {\mathbb N}$ and $k \in {\mathbb Z}_+$ define the space 
$$
{\mathcal H}_k(\varphi_{\nu})= \{f \in H({\mathbb C}^n): 
$$
$$
{\cal N}_{\varphi_{\nu, k}}(f) = 
\sup_{z \in {\mathbb C}^n} 
\frac 
{\vert f(z)\vert (1 + \Vert z \Vert)^k}
{e^{(\psi_{\nu}^*)^*(\ln (1 + \vert Im z_1 \vert), \ldots , \ln (1 + \vert Im z_n \vert ))}} < \infty \}. 
$$
Let for each $\nu \in {\mathbb N}$ \ ${\mathcal H}(\varphi_{\nu})$ be a projective limit of spaces 
${\mathcal H}_k(\varphi_{\nu})$. 
Define the space ${\mathcal H}(\varPhi)$ as an inductive limit of spaces ${\mathcal H}(\varphi_{\nu})$. 

The following assertions allow to obtain an additional information on a structure of the space $E(\varPhi)$.

\begin{proposition}
Let functions of the family $\varPhi$ satisfy the condition $i_2)$ of Theorem 2 and functions $\psi_{\nu}$ be convex on ${\mathbb R}^n$ ($\nu \in {\mathbb N}$). 
Then ${\mathcal H}(\varPhi) = E(\varPhi)$. 
\end{proposition}

\begin{proposition}
Let the family $\varPhi$ satisfies the conditions of Theorem 2. 
Then $E(\varPhi) = {\mathcal H}(\varPhi)$.
\end{proposition}

Next, for each $\nu \in {\mathbb N}$ and $m \in {\mathbb Z}_+$ define the normed space
$$
G_m({\psi_{\nu}^*}) = \{f \in C^m({\mathbb R}^n): 
\Vert f \Vert_{m, \psi_{\nu}^*}
 = \sup_{x \in {\mathbb R}^n, \vert \alpha \vert \le m, \beta \in {\mathbb Z}_+^n}  
\frac 
{\vert x^{\beta}(D^{\alpha}f)(x) \vert}
{\beta! e^{-\psi_{\nu}^*(\beta)}} < \infty \}.
$$
Let $G({\psi_{\nu}^*})= \bigcap \limits_{m = 0}^{\infty} G_m({\psi_{\nu}^*})$,
$G({\Psi^*})= \bigcup \limits_{\nu = 1}^{\infty} G({\psi_{\nu}^*})$. 
Endow $G({\psi_{\nu}^*})$ with a topology defined by the family of norms 
$\Vert \cdot \Vert_{m, \psi_{\nu}^*}$ ($m \in {\mathbb Z}_+$).
Supply $G({\Psi^*})$ with a topology of an inductive limit of spaces $G({\psi_{\nu}^*})$.

For $f \in E(\varPhi)$ define its Fourier transform $\hat f$ by the formula
$$
\hat f(x) = \int_{{\mathbb R}^n} f(\xi) e^{-i \langle x, \xi \rangle} \ d \xi , \ x \in {\mathbb R}^n.
$$

\begin{theorem} 
Let the family $\varPhi$ satisfies  the condition $i_2)$ of Theorem 2 and the following two conditions:

$i_3)$. for each $m \in {\mathbb N}$ there is a constant $a_m > 0$ such that 
$$
\varphi_m(2x) \le \varphi_{m+1}(x) + a_m, \ x \in {\mathbb R}^n;
$$

$i_4)$. for each $m \in {\mathbb N}$ there is a constant $l_m > 0$ such that 
$$
2\varphi_m(x) \le \varphi_{m + 1}(x) + l_m, \ x \in {\mathbb R}^n.
$$
 
Then Fourier transformation ${\mathcal F}: f \in E(\varPhi) \to  \hat f$ establishes an isomorphism of spaces $E(\varPhi)$ and $G({\Psi^*})$.
\end{theorem}

Further, let ${\varPhi^*}=\{{\varphi_{\nu}^*}\}_{\nu=1}^{\infty}$. 
For each $\nu \in {\mathbb N}$ and $m \in {\mathbb Z}_+$ define the space
$$
GS_m(\varphi_{\nu}^*) = \{f \in C^m({\mathbb R}^n): q_{m, \nu}(f) = 
\sup_{x =(z_1, \ldots , x_n) \in {\mathbb R}^n, \atop \vert \alpha \vert \le m} 
\frac {\vert (D^{\alpha}f)(x) \vert} {e^{-\varphi_{\nu}^*(\vert x_1 \vert, \ldots , \vert x_n \vert)}} < \infty \}.
$$
For each $\nu \in {\mathbb N}$ let
$GS({\varphi_{\nu}^*}) = \bigcap \limits_{m \in {\mathbb Z_+}}GS_m({\varphi_{\nu}^*})$. 
Endow $GS({\varphi_{\nu}^*})$ with a topology defined by the family of norms $q_{\nu, m}$ ($m \in {\mathbb Z}_+$). 
Let 
$GS({\varPhi^*})= \bigcup \limits_{\nu \in {\mathbb N}}GS({\varphi_{\nu}^*})$. Supply $GS({\varPhi^*})$ with an inductive limit topology of spaces $GS({\varphi_{\nu}^*})$. 

\begin{theorem}
Let functions of the family $\varPhi$ be convex and satisfy the condition $i_2)$ of Theorem 2 and the condition 
$i_3)$ of Theorem 3. Then $G({\Psi^*}) = GS({\varPhi^*})$.
\end{theorem}

{\bf Acknowledgements}. 
The research was supported by grants from RFBR (14-01-00720, 14-01-97037) and RAS Program of Fundamental Research "`Modern problems of theoretical mathematics"'.


\end{document}